\documentclass{amsart}
\usepackage{amsmath}
\usepackage{amsfonts}
\usepackage{amsthm, upref}
\usepackage{graphicx}
\usepackage{enumerate}
\usepackage[usenames, dvipsnames]{color}
\usepackage{mathtools}
\usepackage{float}

\input xy
\xyoption{all}

\usepackage{ifthen}
\usepackage{tikz}
\usepackage{appendix}
\usepackage{verbatim}
\usetikzlibrary{decorations.pathmorphing}
\usetikzlibrary{calc}
\usepackage{hyperref}
\usepackage{mathtools}

\renewcommand{\comment}[1]{}
\newcommand{\eq}{\begin{equation}}
\newcommand{\en}{\end{equation}}
\newcommand{\rr}{\mathbb{R}}
\newcommand{\NN}{\mathbb{N}}

\newcommand{\abs}[1]{\left\lvert #1 \right\rvert}

\newcommand{\combi}[2]{\begin{pmatrix}#1 \\ #2 \end{pmatrix}}

\newcommand{\E}{\mathrm{E}}

\newcommand{\Ent}{\mathrm{Ent}}

\newcommand{\perm}{\mathbb{S}}
\newcommand{\Hil}{\mathcal{H}}
\newcommand{\trho}{\tilde\rho}
\newcommand{\Fdet}{\mathrm{det}_{\mathrm{F}}}

\begin{document}

\theoremstyle{plain}
\newtheorem{thm}{Theorem}
\newtheorem{lemma}[thm]{Lemma}
\newtheorem{prop}[thm]{Proposition}
\newtheorem{cor}[thm]{Corollary}

\theoremstyle{definition}
\newtheorem{defn}{Definition}
\newtheorem{asmp}{Assumption}
\newtheorem{notn}{Notation}
\newtheorem{prb}{Problem}

\theoremstyle{remark}
\newtheorem{rmk}{Remark}
\newtheorem{exm}{Example}
\newtheorem{clm}{Claim}

\title[Partition function for the Mallows model]{Limiting partition function for the Mallows model: a conjecture and partial evidence}
\author{Soumik Pal}
\address{Department of Mathematics\\ University of Washington\\ Seattle, WA 98195}
\email{soumik@uw.edu}



\thanks{This research is partially supported by NSF grant DMS-2052239, DMS-2134012 and the PIMS PRN-$01$ Kantorovich Initiative.}

\date{\today}

\begin{abstract} Let $\perm_n$ denote the set of  permutations of $n$ labels. We consider a class of Gibbs probability models on $\perm_n$ that is a subfamily of the so-called Mallows model of random permutations. The Gibbs energy is given by a class of right invariant divergences on $\perm_n$ that includes common choices such as the Spearman foot rule and the Spearman rank correlation. Mukherjee \cite{Mukherjee16} computed the limit of the (scaled) log partition function (i.e. normalizing factor) of such models as $n\rightarrow \infty$. Our objective is to compute the exact limit, as $n\rightarrow \infty$, without the log. We conjecture that this limit is given by the Fredholm determinant of an integral operator related to the so-called Schr\"odinger bridge probability distributions from optimal transport theory. We provide partial evidence for this conjecture, although the argument lacks a final error bound that is needed for it to become a complete proof. 
 \end{abstract}

\maketitle

\section{Introduction}
Let $c:[0,1]^2 \rightarrow [0, \infty)$ denote a \textit{cost} function satisfying the following assumptions
\begin{itemize}
\item $c$ is twice continuously differentiable on $[0,1]^2$.
\item $c(x,x)=0$ for all $x\in [0,1]$.
\item $c$ is symmetric, i.e. $c(x,y)=c(y,x)$ for all $(x,y)\in [0,1]^2$.
\item $c(x,y)=c(1-x, 1-y)$ for all $(x,y)\in [0,1]^2$.
\end{itemize}
An example of such a cost function is $c(x,y)=(x-y)^2$. The first of these assumptions is made for technical convenience as will be apparent below. No attempt has been made to get the optimal set of assumptions.

Fix $n\in \NN$. Let $\perm_n$ denote the set of all permutations of $n$ labels $[n]:=\{1,2,\ldots, n\}$. Consider the following quantity 
\begin{equation}\label{eq:paritionfnln}
L_n = \frac{1}{n!} \sum_{\sigma \in \perm_n} \exp\left( - \sum_{i=1}^n c(i/n, \sigma_i/n) \right). 
\end{equation}
We are interested in the limit of this sequence as $n\rightarrow \infty$. The reason it comes up is that this is the partition function of a family of probability distributions which is a subset of the well-known Mallows models \cite{Mallows57} of random permutations. See the Introduction in \cite{Mukherjee16} and many applications listed in \cite[Chapters 5 and 6]{Diaconisbook}. For example, the case of $c(x,y)=(x-y)^2$ is related to the Spearman rank correlation. 

Our goal in this paper is to understand $\lim_{n\rightarrow \infty} L_n$. This problem is important in statistical estimation \cite{Mukherjee16} and also to understand scaling limits of large random permutations with fixed patterns \cite{KKRW20}. See also \cite[Section 2e]{DiaconisRam} for generalizations to other groups where the importance of this problem is stressed. We will try to convince the reader that there are constants $\Gamma_0$ and $C$ such that 
\begin{equation}\label{eq:paritionlimit}
\lim_{n\rightarrow \infty} e^{n \Gamma_0} L_n =C.
\end{equation}  
The value of the constant $\Gamma_0$ is already known due to \cite{Mukherjee16} and is the value of an entropy-regularized optimal transport problem with uniform marginals. See also \cite{Starr09} for a special case of a discontinuous cost function. We conjecture in this paper, and give partial evidence, that, under suitable assumptions, the constant $C$ is the Fredholm determinant  \cite[Definition 14.35]{van22} of a certain integral operator related to the so-called Schr\"{o}dinger bridge, the optimal coupling for the same entropy-regularized optimal transport problem. Both these concepts are described below. Taken together, they give the limiting partition function of this class of Mallows models that satisfy all our assumptions.

The concept of entropy-regularized optimal transport and the related notion of Schr\"odinger bridges can be found in \cite{L12}.  Let $\mu$ denote the Uni$(0,1)$ distribution. Let $\Pi(\mu, \mu)$ denote the set of couplings (i.e., joint distributions) with both marginals $\mu$. Then the entropic OT problem is given as the solution to the following optimization problem on $\Pi(\mu, \mu)$: 
\begin{equation}\label{eq:entropicot}
\Gamma_0:=\inf_{\xi \in \Pi(\mu, \mu)} \left[ \int c(x,y) \xi(x,y) dxdy + \Ent(\xi) \right], 
\end{equation}
where $\Ent(\cdot)$ is the optimal transport entropy (the negative of the usual differential Shannon entropy) given by $\Ent(\xi)= \int \xi(x,y) \log \xi(x,y) dxdy$ if $\xi$ has a density (also denoted by $\xi$) and infinity otherwise. 

The optimal $\rho\in \Pi(\mu, \mu)$ that attains $\Gamma_0$ exists and is called the (static) Schr\"odinger bridge for the cost $c$ and marginals $\mu$ and $\mu$. From the work of R\"uschendorf and Thomsen \cite{RT93} (building on Csiszar \cite{C75}) it is known that the Schr\"odinger bridge always admits a density is of the following form 
\begin{equation}\label{eq:schrodingerform}
\rho(x,y) = \exp\left( - c(x,y) - a(x) - a(y) \right),
\end{equation}
for some measurable function $a$ satisfying the following marginal constraint almost surely.
\begin{equation}\label{eq:marginalconstraint}
 \int_0^1 e^{- c(x,y) - a(y)}  dy=e^{-a(x)}, \quad \text{for}\;  x\in [0,1]. 
\end{equation}
Please note that we are using a standard abuse of the notation by referring to both the measure and its density by the letter $\rho$.

In particular, $\rho$ is symmetric in its argument (since both $c$ and the marginal constraints are symmetric in the coordinates) and 
\[
\begin{split}
\Gamma_0 &=  \int c(x,y) \rho(x,y) dxdy + \Ent(\rho) \\
&= \int c(x,y) \rho(x,y) dxdy - \int \left( c(x,y) + a(x) + a(y) \right) \rho(x,y)dxdy\\
&= - 2 \int_0^1 a(x)dx,
\end{split}
\]
where the final equality is due to the fact that $\rho \in \Pi(\mu, \mu)$. 

Mukherjee \cite[Theorem 1.5]{Mukherjee16} shows that the log-partition function has the following large deviation limit
\[
\lim_{n\rightarrow \infty}\frac{1}{n} \log L_n = -\Gamma_0. 
\]
To compare our notation with that of \cite{Mukherjee16}, note that $\theta=1$, their $f=-c$, $\mathcal{M}=\Pi(\mu, \mu)$, $D(\cdot || u)= \Ent(\cdot)$, $Z_n(f, \theta)=\log (n! L_n)$ and $Z_n(0)=\log n!$.

Hence, it makes sense to consider $\lim_{n\rightarrow \infty} e^{n \Gamma_0} L_n$. Towards that goal, define 
\begin{equation}\label{eq:denom1}
\begin{split}
D_n &:= \frac{1}{n!} \sum_{\sigma \in \perm_n} \prod_{i=1}^n \rho(i/n, \sigma_i/n) \\
&=\frac{1}{n!} \sum_{\sigma \in \perm_n} \exp\left( - \sum_{i=1}^n c(i/n, \sigma_i/n) - 2\sum_{i=1}^n a(i/n) \right)=L_n \exp\left( - 2 \sum_{i=1}^n a(i/n) \right)\\
&\approx L_n \exp\left( - 2n \int_0^1 a(x) dx\right)=e^{n \Gamma_0} L_n.
\end{split}
\end{equation}
The $\approx$ in the middle can be quantified as the discretization error in the Riemann sum approximation. In fact, assuming that $a$ is twice continuously differentiable, we get 
\[
\begin{split}
\int_0^1 a(x)dx &- \frac{1}{n}\sum_{i=1}^n a(i/n) = \sum_{i=1}^n \int_{(i-1)/n}^{i/n}  \left( a(x) - a(i/n) \right) dx \\
&= \sum_{i=1}^n a'(i/n) \int_{(i-1)/n}^{i/n} (x- i/n)dx + O\left( \sum_{i=1}^n \int_{(i-1)/n}^{i/n} (x-i/n)^2 dx \right)\\
&=  \frac{1}{2n^2} \sum_{i=1}^n a'(i/n) + O\left( \frac{1}{n^2} \right).
\end{split}
\]
How can we guarantee that $a$ is twice continuously differentiable? This follows from the assumed twice  continuous differentiability of $c$ and the integral equation \eqref{eq:marginalconstraint}. Hence, 
\[
\lim_{n\rightarrow \infty} n\left[ \int_0^1 a(x)dx - \frac{1}{n}\sum_{i=1}^n a(i/n) \right] = \frac{1}{2}\int_0^1 a'(x)dx = \frac{a(1)-a(0)}{2}. 
\]
Now, due to our assumption $c(x,y)=c(1-x, 1-y)$, we must have $a(0)=a(1)$. This is because the invariance of $\mu$ under the map $x\mapsto 1-x$. Hence, if $\rho(x,y)$ is the Schr\"odinger bridge, so is $\rho(1-x, 1-y)$ due to the uniqueness of the solution of the strictly convex optimization problem \eqref{eq:entropicot}.

Hence, from \eqref{eq:denom1}, 
\begin{equation}\label{eq:denom2}
\lim_{n\rightarrow \infty} e^{n \Gamma_0} L_n= \lim_{n\rightarrow \infty} D_n.
\end{equation}
This is what we will evaluate below. 

The first step is to consider the Markov integral operator corresponding to the probability density $\rho$ and derive its spectral decomposition. Consider the separable Hilbert space $\Hil=L^2[0,1]$. Define the integral operator: for $u \in \Hil$, 
\[
Tu(x):= \int_0^1 u(y) \rho(x,y)dy.  
\]
Clearly, for any other $v\in \Hil$, by Fubini's Theorem and the symmetry of $\rho$, 
\[
\begin{split}
\int_0^1 v(x) Tu(x) dx &= \int_0^1 \int_0^1 v(x) \rho(x,y) u(y) dy dx= \int_0^1 u(y) \left[ \int_0^1 v(x) \rho(y,x) dx \right] dy\\
&= \int_0^1 u(y) Tv(y) dy. 
\end{split}
\]
In particular, $T$ is a self-adjoint linear operator on the separable Hilbert space $\Hil$. Since this is a Hilbert-Schmidt operator, it is also compact. Hence, it admits a spectral decomposition. In particular, there exists a countable sequence of eigenvalues $(\lambda_n,\; n \in \NN )$ and their corresponding eigenfunctions (with multiplicities) such that 
\[
\rho(x,y) = 1 + \sum_{n=1}^ \infty \lambda_n \phi_n(x)\phi_n(y). 
\]
All the eigenvalues are real since $T$ is self-adjoint. In fact, $I - T$ is a nonnegative operator, where $I$ is the identity operator. Thus all the eigenvalues of $T$ lie in the interval $[-1,1]$. The eigenfunction corresponding to eigenvalue $1$ is the constant function $\phi_1(x) \equiv 1$.

\begin{asmp}\label{asmp:spectralgap}
Assume that there is a positive spectral gap, i.e., 
\begin{equation}\label{eq:spectralgap}
\sigma^2 = \left(1- \max_{i\ge 1} \lambda^2_i\right) \in (0,1). 
\end{equation}
\end{asmp}

\begin{asmp}\label{asmp:lipefun}
Assume that each eigenfunction $\phi_n$ is Lipschitz continuous on $[0,1]$.
\end{asmp}

Consider $T$ as a self-adjoint operator on $\Hil_1$, the subspace of $\Hil$ that is orthogonal to the constant functions. Since $T^2$ is a trace class operator on the Hilbert space $\Hil_1$, the Fredholm determinant of $I - T^2$ exists and is given by the absolutely convergent infinite product 
\[
\Fdet(I - T^2)= \prod_{n=1}^\infty(1- \lambda^2_n). 
\]
See \cite[Definition 14.35, Theorem 14.44]{van22}.
\medskip

\noindent\textbf{Conjecture.} Our main conjecture is that, under Assumptions \ref{asmp:spectralgap} and \ref{asmp:lipefun}, the following limit holds.
\begin{equation}\label{eq:conjecture}
\lim_{n\rightarrow \infty} D_n = C= \frac{1}{\sqrt{\Fdet(I - T^2)}}. 
\end{equation}
Hence, by \eqref{eq:denom2}, $\lim_{n\rightarrow \infty} e^{n \Gamma_0} L_n= \left(\Fdet(I - T^2)\right)^{-1/2}$ which completes our aim outlined in \eqref{eq:paritionlimit}. 

Although we will give a partial proof towards this conjecture in the next section, let us provide some intuition why such a limit should be true. In \cite[Theorem 2]{HLP}, the present author and coauthors proved a similar but more complex limit. The relationship between this paper and that one may be explained in the following way. If we consider $L_n$ in \eqref{eq:paritionfnln} as a function of the empirical distribution $\hat{\mu}_n:=\frac{1}{n} \sum_{i=1}^n \delta_{i/n}$, in \cite{HLP} we consider a similar function of the empirical distribution $\tilde{\mu}_n := \frac{1}{n} \sum_{i=1}^n \delta_{X_i}$, where $X_1, X_2, \ldots$ is a sequence of i.i.d. Uni$(0,1)$ random variables. In the latter case, the limit $\lim_{n\rightarrow\infty} D_n$ is random and belongs to the class of second order Gaussian chaos as shown in \cite[Theorem 2]{HLP}. The reason we get the Gaussian chaos is due to the limiting Gaussian fluctuation between $\tilde{\mu}_n$ and Uni$(0,1)$ as established by standard empirical process theory.  There is, of course, no limiting Gaussian fluctuation for the difference between $\hat{\mu}_n$ and Uni$(0,1)$. Hence one may assume that the limiting Gaussian random variables all have zero variance. If we plug this in \cite[Theorem 2]{HLP} and simplify to our case at hand we get \eqref{eq:conjecture}.

Although this connection has been pointed out in the introduction of \cite{HLP}, that proof simply cannot cover this case due to the lack of randomness. The difference between the two set-ups may be explained by the following analogy. Whereas the proof in \cite{HLP} can generalize to sampling with replacement from the finite set $(i/n,\; i \in [n])$, our current set-up is about sampling without replacement. The combinatorics is much more involved which leads to our inability to completing the proof of the conjecture.  
\bigskip

\noindent\textbf{Addendum in Oct 2026.} This article was posted on the arxiv in June 2024. I am happy to report that this conjectured has now been settled. In fact, there have been three independent groups that have successfully settled this problem. The first proof appeared in  \cite{tripathi26} (May '26) and is due to Raghavendra Tripathi. Neither the condition $c(x,x)\equiv 0$ nor $c(x,y)=c(1-x,1-y)$ are needed for the result to hold. The second proof appeared shortly afterwards (also May '26) and is due to Ottolini and Starr \cite{OS26}. This paper shows that the conjecture remains true even for discontinuous cost functions that are constant over blocks. The symmetry condition $c(x,y)=c(y,x)$ is also dropped. The third proof appears in the work by Duan, Lyu and Powell \cite{DLP26} (October '26) where they show that the conjecture remains true for $C^1$ costs (see Corollary 2.15). More importantly, they uncover generalizations of this conjecture to Poisson contingency tables.

\section{A partial proof of the conjectured limit} 

Let $\trho$ be the kernel $\rho-1$. Then $\trho(x,y)= \sum_{i=1}^\infty \lambda_i \phi_i(x)\phi_i(y)$ where the series converges in $L^2$. In particular, due to the marginal constraints, 
\begin{equation}\label{eqn:marginal}
    \int_0^1 \tilde{\rho}(x,y)dy=0=\int_0^1 \tilde{\rho}(z,w)dz, 
\end{equation}
for $x,w$ in $[0,1]$. 

For any choice of $(x_1, \ldots, x_n)$ and any $\sigma \in \perm_n$,
\[
\begin{split}
\prod_{i=1}^n \rho(x_i, x_{\sigma_i}) &= \prod_{i=1}^n \left(1 + \trho(x_i, x_{\sigma_i}) \right)= 1 + \sum_{A \subseteq [n],\; A\neq \emptyset} \prod_{i\in A} \trho(x_i, x_{\sigma_i}).
\end{split}
\]
Hence, 
\[
\begin{split}
D_n &= \frac{1}{n!} \sum_{\sigma \in \perm_n} \left[1 + \sum_{A \subseteq [n],\; A\neq \emptyset} \prod_{i\in A} \trho(i/n, \sigma_i/n)\right]\\
&= 1 + \frac{1}{n!} \sum_{\sigma \in \perm_n}  \sum_{A \subseteq [n],\; A\neq \emptyset} \prod_{i\in A} \trho(i/n, \sigma_i/n)\\
&= 1 + \frac{1}{n!} \sum_{\sigma \in \perm_n} \sum_{r=1}^n \sum_{A: \abs{A}=r} \prod_{i \in A} \trho(i/n, \sigma_i/n)\\
&= 1 + \frac{1}{n!} \sum_{r=1}^n \sum_{A: \abs{A}=r} \sum_{\sigma \in \perm_n} \prod_{i \in A} \trho(i/n, \sigma_i/n)\\
&= 1 + \frac{(n-r)!}{n!} \sum_{r=1}^n \sum_{1\le i_1 < i_2 < \cdots < i_r \le n} \sum_{1\le j_1 \neq j_2 \neq \cdots \neq j_r \le n} \prod_{t=1}^r \trho(i_t/n, j_t/n).
\end{split}
\]
Here, the condition $\{1\le j_1 \neq j_2 \neq \cdots \neq j_r \le n\}$ means all the indices are distinct and in $[n]$.

Fix $K \in \NN$. For all $n\ge K$, let 
\[
D_{n,K}:= 1 + \frac{(n-r)!}{n!} \sum_{r=1}^K \sum_{1\le i_1 < i_2 < \ldots i_r \le n} \sum_{1\le j_1 \neq j_2 \neq \ldots \neq j_r \le n} \prod_{t=1}^r \trho(i_t/n, j_t/n).
\]
For $L\in \NN$, define 
\[
\trho^{(L)}(x,y) = \sum_{l=1}^L \lambda_l \phi_l(x) \phi_l(y).
\]

Finally, define 
\[
\begin{split}
D_{n,K}^{(L)}&:= 1 + \frac{(n-r)!}{n!} \sum_{r=1}^K \sum_{1\le i_1 < i_2 < \cdots < i_r \le n} \sum_{1\le j_1 \neq j_2 \neq \cdots \neq j_r \le n} \prod_{t=1}^r \trho^{(L)}(i_t/n, j_t/n)\\
&= 1 + \frac{(n-r)!}{n!} \sum_{r=1}^K \sum_{1\le i_1 < i_2 < \cdots < i_r \le n} \sum_{1\le j_1 \neq j_2 \neq \cdots \neq j_r \le n} \prod_{t=1}^r \left[ \sum_{l=1}^L \lambda_l \phi_l(i_t/n) \phi_l(j_t/n)\right]\\
&= 1 + \sum_{r=1}^K \sum_{1\le l_1, \ldots, l_r \le L} \prod_{t=1}^r \lambda_{l_t} \times \\
&\left[ \frac{(n-r)!}{r! n!} \sum_{1\le i_1 \neq \cdots \neq i_r \le n} \sum_{1\le j_1 \neq \cdots \neq j_r \le n} \prod_{t=1}^r \phi_{l_t}(i_t/n) \phi_{l_t}(j_t/n)\right]. 
\end{split}
\]
Recall that the sum over indices $i_1\neq \cdots \neq i_r$ means that all the indices are distinct. 

Now fix a vector $(l_1, l_2, \ldots, l_r)$ in $[L]^r$. Let $a_l$ denote the frequency of appearance of $l\in [L]$ in this sequence. Then each $a_l\ge 0$ and assume $\sum_{l=1}^L a_l =r\le K$. Since, for fixed $r$, 
\[
\frac{(n-r)!}{n!} \approx \frac{1}{n^r}, 
\]
consider the normalized inner sum
\begin{equation}\label{eq:perfectsq}
\begin{split}
    \frac{1}{n^r} &\sum_{1\le i_1 \neq \cdots \neq i_r \le n} \sum_{1\le j_1 \neq \cdots \neq j_r \le n} \prod_{t=1}^r \phi_{l_t}(i_t/n) \phi_{l_t}(j_t/n)\\
    &=\left[ \frac{1}{n^{r/2}} \sum_{1\le i_1 \neq \cdots \neq i_r \le n}  \prod_{t=1}^r \phi_{l_t}(i_t/n) \right]^2
\end{split}
\end{equation}
The claim is unless $a_l\in \{0, 2\}$ for all $l\in [L]$, the contributions of the corresponding terms in the normalized sum inside the square converge to zero as $n\rightarrow \infty$. 

To see this let us introduce a sequence of i.i.d. Uni$(0,1)$ random variables $(U_1, U_2, \ldots)$. For $n \in \NN$, define 
\[
U_i^{(n)}:=\frac{1}{n} \lceil n U_i \rceil.
\]
Then, for each $n$, $\left( U_1^{(n)}, U^{(2)}_n, \ldots \right)$ is a sequence of i.i.d. discrete Uni$[n]$ random variables, and obviously, $\abs{U_i^{(n)} - U_i}\le 1/n$. 

Note  
\[
\frac{1}{n^r}\sum_{1\le i_1 \neq \cdots \neq i_r \le n}  \prod_{t=1}^r \phi_{l_t}(i_t/n)=\E\left[\left(\prod_{t=1}^r \phi_{l_t}\left(U^{(n)}_t\right)\right) 1\left\{ U^{(n)}_1 \neq U^{(n)}_2 \neq \cdots \neq U^{(n)}_r\right\} \right].
\]
Thus, to analyze the limit of \eqref{eq:perfectsq}, it suffices to analyze $n^{r/2}$ times the RHS. 

Consider the following events $F^{(n)}_{ij}=\{U^{(n)}_i = U^{(n)}_j\}$, for $i< j \in \NN$. Then 
\[
\cup_{1\le i< j\le r} F^{(n)}_{ij}= \{U^{(n)}_1 \neq U^{(n)}_2 \neq \cdots \neq U^{(n)}_r\}^c.
\]
Thus
\[
\begin{split}
\E&\left[\left(\prod_{t=1}^r \phi_{l_t}(U^{(n)}_t)\right) 1\{ U^{(n)}_1 \neq U^{(n)}_2 \neq \cdots \neq U^{(n)}_r\} \right]\\
=& \E\left[\left(\prod_{t=1}^r \phi_{l_t}(U^{(n)}_t)\right)\right]
-\E\left[\left(\prod_{t=1}^r \phi_{l_t}(U^{(n)}_t)\right); \cup_{1\le i < j \le r} F^{(n)}_{ij} \right].
\end{split}
\]
Here, for any event $A$ and any integrable random variable $Y$ in a probability space, $\E(Y;A):=\E(Y1_A)$.


Since every $\phi_i$ is an eigenfunction orthogonal to $1$, they satisfy $\E(\phi_i(U_1))=0$ and $\E(\phi_i^2(U_1))=1$. 
Due to the spatial discreteness this is not going to be exactly true for $U_i^{(n)}$. However, by Assumption \ref{asmp:lipefun} on the Lipschitzness of eigenfunctions, it follows that 
\begin{equation}\label{eq:exactexpect}
\begin{split}
\E\left( \phi_i(U^{(n)}_1)\right) &= O\left( \frac{1}{n} \right), \\
 \E\left[\left(\prod_{t=1}^r \phi_{l_t}(U^{(n)}_t)\right) \right] &= \prod_{t=1}^r \E\left( \phi_{l_t}(U^{(n)}_t) \right) =O\left( \frac{1}{n^r}\right),
\end{split}
\end{equation} 
where the constant in $O(1/n)$ can be chosen uniformly for any finite collection of eigenfunctions. In the calculation below, every time we encounter an expression as in \eqref{eq:exactexpect} we will ignore the $O(1/n)$ term and put zero instead. This is simply for the clarity of the combinatorial expressions. Because $L$ and $K$ are both finite, only finitely many eigenfunctions ever get used and the constant in $O(1/n)$ remains uniformly bounded. The primary reason why we failed to complete this proof is because we cannot suitably estimate this error when $L$ and $K$ are not bounded. Nevertheless, with this convenient abuse  of notation,
\[
\begin{split}
n^{r/2}\E&\left[\left(\prod_{t=1}^r \phi_{l_t}(U^{(n)}_t)\right) 1\{ U^{(n)}_1 \neq U^{(n)}_2 \neq \cdots \neq U^{(n)}_r\} \right]\\
&=- n^{r/2}\E\left[\left(\prod_{t=1}^r \phi_{l_t}(U^{(n)}_t)\right); \cup_{1\le i < j \le r} F^{(n)}_{ij} \right],
\end{split}
\]
where, let us repeat again, we have ignored an $O(n^{-r/2})$ error.

On the other hand, by the inclusion-exclusion principle, and by utilizing exchangeability,
\[
\begin{split}
    \E&\left[\left(\prod_{t=1}^r \phi_{l_t}(U^{(n)}_t)\right); \cup_{1\le i < j \le r} F^{(n)}_{ij} \right]\\
    &=\sum_{k=1}^{r(r-1)/2} (-1)^{k-1} \sum \E\left[\left(\prod_{t=1}^r \phi_{l_t}(U^{(n)}_t)\right); \text{intersection of $k$ many $F^{(n)}_{ij}$s} \right], 
\end{split}
\]
where the inner sum is over all choices of $k$ many $F^{(n)}_{ij}$s. Fix $(i_1, j_1), \ldots, (i_k, j_k)$. Then 
\[
\cap_{m=1}^k F^{(n)}_{i_m j_m}= \{ U^{(n)}_{i_1}=U^{(n)}_{j_1}, \ldots, U^{(n)}_{i_k}=U^{(n)}_{j_k}\}. 
\]

We now make the following observations. The above constraint gives a partition of $[r]$, where, for each block of the partition, the discrete uniform random variables corresponding to the indices in that block take the same value. 

\begin{itemize}
    \item If the partition contains a block that is a singleton, i.e. $\cup_{m=1}^k \{i_m, j_m\} \neq [r]$, there will be some $U^{(n)}_t$ which has no constraint and is, therefore, independent of the other uniform random variables. Let the number of singletons be $r-r'$, for some $r' \in \{0\} \cup [r-1]$. Then, by \eqref{eq:exactexpect}, independence and exchangeability,  
    \[
    \begin{split}
    \E&\left[\left(\prod_{t=1}^r \phi_{l_t}(U^{(n)}_t)\right);  \cap_{m=1}^k F^{(n)}_{i_mj_m} \right]\\
    =&O\left(\frac{1}{n^{r-r'}}\right)  \E\left[\left(\prod_{t=1}^{r'} \phi_{l_t}(U^{(n)}_t)\right);  \cap_{m=1}^k F^{(n)}_{i_mj_m} \right].
    \end{split}
    \]
    The cases below will show that $n^{r'/2}$ times the expectation on the RHS remains bounded, as $n\rightarrow \infty$. Hence $n^{r/2}$ times the expectation on the LHS goes to zero as $n\rightarrow \infty$ whenever $r' < r$.  
    Hence, asymptotically, the only non-zero terms come from partitions of $[r]$ that do not contain any singleton blocks. 
    \item Now consider the case where every block in the partition is of size two. Such partitions are in correspondence with perfect matchings of the complete graph $K_r$. In particular, $r=2k$ must be even. For any such perfect matching, say $U^{(n)}_{i_1}=U^{(n)}_{j_1}$, $U^{(n)}_{i_2}=U^{(n)}_{j_2}$, etc., 
    \[
    \E\left[\left(\prod_{t=1}^r \phi_{l_t}(U^{(n)}_t)\right);  \cap_{m=1}^k F^{(n)}_{i_mj_m} \right]= \frac{1}{n^{r/2}}\prod_{m=1}^{r/2} \E\left( \phi_{l_{i_m}}(U^{(n)}_{i_m}) \phi_{l_{j_m}}(U^{(n)}_{i_m}) \right).
    \]
    The above product is zero in all cases, except when $l_{i_m}=l_{j_m}$, for all $m\in [r/2]$, in which case the product is one. This is due to the orthonormal property of the eigenfunctions. Both claims hold up to a smaller discretization error as in \eqref{eq:exactexpect}. When multiplied by $n^{r/2}$, each such expectation gives an $O(1)$ term, and there are only finitely many matchings of $K_r$, $r\le K$.  
    
    \item Finally, consider the contribution of a partition of $[r]$ that contains a block of size $3$ or more and no singletons. Then, obviously, $P\left( \cap_{m=1}^k F^{(n)}_{i_mj_m}\right)=o\left( \frac{1}{n^{r/2}}\right)$.
    Since the eigenfunctions are uniformly bounded for $r\le K$, 
    \[
    \E\left[\left(\prod_{t=1}^r \phi_{l_t}(U^{(n)}_t)\right);  \cap_{m=1}^k F^{(n)}_{i_mj_m} \right]=o\left( \frac{1}{n^{r/2}}\right).
    \]
\end{itemize}

Since $r\le K$ and there are only finitely many partitions of $[r]$, 
\[
\sum \E\left[\left(\prod_{t=1}^r \phi_{l_t}(U_t)\right); \text{intersection of $k$ many $F^{(n)}_{ij}$s} \right]=o\left( \frac{1}{n^{r/2}} \right)
\]
for all $r$ odd and all $k\neq r/2$, when $r$ is even. The only remaining case is when $r$ is even and $k=r/2$. This is the computation done in the second bulleted item in the itemized list above. The limiting contribution of each term in that case is either zero or one as shown. Hence the total contribution is the number of terms that contribute one. 

To find this number, let $a_i\in \NN \cup \{0\}$ denote the number of times $l_i$ appears in the sequence $(l_1, \ldots, l_r)$. If any $a_i$ is odd, there is no matching of $K_r$ that matches all $l_i$s to themselves and then the sum is zero. If all $a_i$s are even, the only partitions of $r$ that have nonzero contributions are precisely those that belong to the direct product of the set of perfect matchings of the complete graph $K_{a_i}$. Hence, the number of such terms is the product of the number of perfect matchings of $K_{a_i}$s. Thus, in this only remaining case, 
\[
\sum \E\left[\left(\prod_{t=1}^r \phi_{l_t}(U_t)\right); \text{intersection of $k$ many $F_{ij}$s} \right]=\frac{1}{n^{r/2}} \prod_{i=1}^\infty (a_i-1)!!,
\]
where, by convention, $(0-1)!!=1$. Let $a_i:=2b_i$; the above may also be written as 
\[
\frac{1}{n^{r/2}} \prod_{i=1}^\infty (2b_i-1)!!.
\]

Combining all these terms, 
\begin{equation}\label{eq:denapprox1}
\begin{split}
    D_{n,K}^{(L)}= 1 +  \sum_{r=1,\; \text{even}}^K \frac{1}{r!}\sum \prod_{i=1}^L \lambda_{i}^{2b_i} \left( \prod_{i=1}^L (2 b_i-1)!!\right)^2 + o(1),\;\; \text{as $n\rightarrow \infty$}, 
\end{split}
\end{equation}
where the inner sum is over all sequences $(l_1, \ldots, l_r)\in [L]^r$ and the nonnegative integers $(2b_i,\; i\in [L])$ such that $2\sum_{i=1}^L b_i=r$ record the frequency of appearance of $i$ in the sequence.

We now recall a property of Hermite polynomials, i.e., the number of perfect matchings of $K_n$ is exactly $H_n(0)$ which is the value at zero for the $n$th Hermite polynomial
\[
H_n(x)=(-1)^n e^{x^2/2}\frac{d^n}{dx^x} e^{-x^2/2}.
\]
See \cite[eqn. (3.14)]{HL71} for a more general identity involving moment polynomials. Also see \cite{Godsil81} for similar identities involving more general graphs. 

Thus, one may also write \eqref{eq:denapprox1} as
\begin{equation}\label{eq:denapprox2}
D_{n,K}^{(L)}= 1 +  \sum_{r=1, \; \text{even}}^K  \frac{1}{r!} \sum \left( \prod_{i=1}^L \lambda_{i}^{b_i} H_{2 b_i}(0)\right)^2 + o_n(1),\;\; \text{as $n\rightarrow \infty$},
\end{equation}
where, as before, the inner sum is over all sequences $(l_1, \ldots, l_r)\in [L]^r$ and the nonnegative integers $(2b_i,\; i\in [L])$ such that $2\sum_{i=1}^L b_i=r$ record the frequency of appearance of $i$ in the sequence.

Consider any permutation symmetric function $f:[L]^r\rightarrow \rr$. For every choice of nonnegative integers $\mathbf{a}:=(a_1, \ldots, a_L)$ such that $\sum_{i=1}^L a_i=r$, let $\Gamma(a_1, \ldots, a_L)$ denote all subsets of $[L]^r$ such that $i\in [L]$ appears exactly $a_i$ times. Pick a representation element $\ell_{\mathbf{a}} \in \Gamma(a_1, \ldots, a_L)$. It then follows easily that 
\[
\frac{1}{r!} \sum_{(l_1, \ldots, l_r)\in [L]^r} f(l_1, \ldots, l_r) = \sum_{(a_1, \ldots, a_L)} \frac{1}{a_1 ! \ldots a_L!}f(\ell_{\mathbf{a}}). 
\]
Thus, from \eqref{eq:denapprox2},
\begin{equation}\label{eq:denapprox3}
\lim_{n\rightarrow \infty}D_{n,K}^{(L)}= D_K^{(L)}:= 1 +  \sum_{r=1, \; \text{even}}^K  \sum \frac{1}{(2b_1)! \ldots (2b_L)!} \left( \prod_{i=1}^L \lambda_{i}^{b_i} H_{2 b_i}(0)\right)^2,
\end{equation}
where now the inner sum is over all choices of nonnegative integers $(b_1, \ldots, b_L)$ such that $2\sum_{i=1}^Lb_i=r$.

Note that, at least formally, 
\[
\begin{split}
\lim_{K\rightarrow \infty} D_K^{(L)}&=D^{(L)}:=1 +  \sum \prod_{i=1}^L \frac{ \lambda_{i}^{2b_i}}{(2b_i)!} \left( H_{2 b_i}(0)\right)^2.
\end{split}
\]
where now the sum is over all choices of nonnegative integers $(b_1, \ldots, b_L)$. The fact that the expression on the RHS is finite (and therefore the limit exists) is a consequence of the  multilinear Mehler formula for the Hermite polynomials. See \cite[eqn. (2.4)]{Foata81} for a choice of $n=L$, $S_n$ to be the $L\times L$ zero matrix, the diagonal matrix $D_n$ to have the diagonal vector $(\lambda^2_1, \ldots, \lambda^2_L)$ and the indeterminate vectors $y,z$ to both be the zero vector. For this choice, in their notation, the only symmetric matrices $N$ that will contribute nonzero terms must be diagonal with $\nu_{ii}=a_i$ equal to our $2b_i$. Hence, by \cite[eqn. (2.4)]{Foata81} the following limit exists and is given by the simple determinantal expression
\[
D^{(L)}=\frac{1}{\prod_{i=1}^L\sqrt{1-\lambda^2_i}}.
\]
Note that this is finite by our assumption on the positive spectral gap \eqref{eq:spectralgap}. Therefore, 
\[
\lim_{L\rightarrow \infty} D^{(L)}:= \frac{1}{\prod_{i=1}^\infty\sqrt{1-\lambda^2_i}}=\frac{1}{\sqrt{\Fdet\left(I - T^2\right)}}.
\]
The relation of the above with the limit conjectured in \eqref{eq:conjecture} is now obvious. What we have shown is that 
\[
\lim_{K,L\rightarrow \infty} \lim_{n\rightarrow \infty} D^{(L)}_{n,K}= \frac{1}{\sqrt{\Fdet\left(I - T^2\right)}},
\]
while what we need to show is 
\[
\lim_{n\rightarrow \infty} \lim_{K,L\rightarrow \infty}  D^{(L)}_{n,K}= \frac{1}{\sqrt{\Fdet\left(I - T^2\right)}}.
\]
The interchange of limits requires a uniform error bound. In \cite{HLP} such an error bound has been proved for that set-up. But the argument does not extend to this case because the combinatorics is different. However, it is reasonable to guess that a more careful combinatorics will provide us with the error bound to establish \eqref{eq:conjecture}.

\section{Acknowledgment} I am very grateful to the editors for giving me an opportunity to contribute towards IJPAM volume celebrating the memory of Prof. K. R. Parthasarathy. May his work continue to inspire many further generations of Indian probabilists. I also thank an anonymous referee for a thorough reading and many comments that improved the paper.

\bibliographystyle{amsalpha}
\bibliography{references}

@article{Godsil81, 
    author={Godsil, C. D.},
    title={Hermite polynomials and a duality relationship for matching polynomials},
    journal={Combinatorica},
    volume={1},
    number={3},
    year={1981},
    pages={257--262}
}

@article{DiaconisRam,
author = {Persi Diaconis and Arun Ram},
title = {{Analysis of systematic scan Metropolis algorithms using Iwahori-Hecke algebra techniques.}},
volume = {48},
journal = {Michigan Mathematical Journal},
number = {1},
publisher = {University of Michigan, Department of Mathematics},
pages = {157 -- 190},
year = {2000},
doi = {10.1307/mmj/1030132713},
URL = {https://doi.org/10.1307/mmj/1030132713}
}

@article{Starr09,
    author = {Starr, Shannon},
    title = "{Thermodynamic limit for the Mallows model on $S_n$ }",
    journal = {Journal of Mathematical Physics},
    volume = {50},
    number = {9},
    pages = {095208},
    year = {2009},
    month = {07},
    issn = {0022-2488},
    doi = {10.1063/1.3156746},
    url = {https://doi.org/10.1063/1.3156746},
    eprint = {https://pubs.aip.org/aip/jmp/article-pdf/doi/10.1063/1.3156746/13814524/095208\_1\_online.pdf},
}

@unpublished{tripathi26,
	author={Tripathi, Raghavendra},
	title={On the partition function of a class of {M}allows model},
	year={2026},
	note={Available at ar{X}iv:2605.03647v1 [math.{PR}] }
}

@unpublished{OS26,
	author={Ottolini, Andrea and Starr, Shannon},
	title={Pal's permanent conjecture: proof for block uniform matrices},
	year={2026},
	note={Available at ar{X}iv:2605.25274v1 [math.{CO}] }
}

@unpublished{DLP26,
	author={Duan, Danny and Lyu, Hanbaek and Powell, William},
	title={Poisson contingency tables, {G}ibbs permutations, and entropic optimal transport},
	year={2026},
	note={Available at ar{X}iv:2610.14616v1 [math.{PR}]}
}

@book{van22,
  title={Functional Analysis},
  author={van Neerven, J.},
  isbn={9781009232494},
  series={Cambridge Studies in Advanced Mathematics},
  url={https://books.google.com/books?id=ln92EAAAQBAJ},
  year={2022},
  publisher={Cambridge University Press}
}

@article{Mukherjee16,
	author={Mukherjee, Sumit},
	title={Estimation in exponential families on permutations},
	journal={The Annals of Statistics},
	volume={44},
	number={2},
	year={2016},
	pages={853--875}
}

@article{Mallows57,
	author={Mallows, C. L.},
	title={Non-null ranking models. {I}.},
	journal={Biometrika},
	volume={44},
	pages={114--130},
	year={1957},
}

@article{KKRW20,
	author={Kenyon, R. and Kr\'{a}l', D. and Radin, C. and Winkler, P.},
	title={Permutations with fixed pattern densities},
	journal={Random Structures and Algorithms},
	volume={56},
	year={2020}
}

@book{Diaconisbook, 
	author={Diaconis, P.},
	title={Group Representations in Probability and Statistics},
	year={1988},
	publisher={IMS, Hayward, CA},
	series={Institute of Mathematics Statistics Lecture Notes -- Monograph Series},
	volume={11},
}

@unpublished{HLP, 
	author={Harchaoui, Z. and Liu, L. and Pal, S.},
	title={Asymptotics of discrete {S}chr{\"o}dinger bridges via chaos decomposition},
	year={2020},
	note={To appear in \textit{Bernoulli}. Preprint available at \texttt{arxiv.org/abs/2011.08963}.}
}

@article{RT93,
title={Note on the {S}chr\"odinger equation and {I}-projections},
author={R\"uschendorf, L. and Thomsen, W.},
journal={Statistics and Probability Letters},
volume={17},
pages={369--375},
year={1993}
}

@article{C75,
title={I-divergence geometry of probability distributions and minimization problems},
author={Csiszar, I.},
journal={Ann. Probab.},
volume={3},
pages={146--158},
year={1975}
}

@article{L12,
  title={From the {S}chr{\"o}dinger problem to the {M}onge-{K}antorovich problem},
  author={L{\'e}onard, Christian},
  journal={Journal of Functional Analysis},
  volume={262},
  number={4},
  pages={1879--1920},
  year={2012},
  publisher={Elsevier}
}

@article{HL71,
    author={Heilmann, Ole J. and Lieb, Elliott H.},
    title={Theory of Monomer-Dimer systems},
    journal={Commun. math. Phys.},
    volume={25},
    pages={190--232},
    year={1972}
}

@article{Foata81,
title = {Some {H}ermite polynomial identities and their combinatorics},
journal = {Advances in Applied Mathematics},
volume = {2},
number = {3},
pages = {250-259},
year = {1981},
issn = {0196-8858},
doi = {https://doi.org/10.1016/0196-8858(81)90006-3},
url = {https://www.sciencedirect.com/science/article/pii/0196885881900063},
author = {Dominique Foata}
}

\end{document}